\documentclass{article}

\usepackage{spconf}
\usepackage{cite}
\usepackage{amsmath,amssymb,amsthm}
\usepackage{mathtools}
\usepackage{bm}
\usepackage{hyperref}

\newtheorem*{theorem*}{Theorem}
\newtheorem{definition}{Definition}
\newtheorem*{proof*}{Proof}

\title{On the Equivalence of \MakeLowercase{$f$}-Divergence Balls and Density Bands \\ in Robust Detection}
\name{Michael Fau\ss{}\textsuperscript{$\ast$}, Abdelhak M.~Zoubir\textsuperscript{$\ast$}, and H.~Vincent Poor\textsuperscript{$\dag$}\thanks{This work was supported in part by the U.S.~National Science Foundation under Grants CCF-1420575  and ECCS-1549881.}}
\address{\textsuperscript{$\ast$}Signal Processing Group, Technische Universit\"at Darmstadt, 64283 Darmstadt, Germany \\[0.2ex]
         \textsuperscript{$\dag$}Dept.~of Electrical Engineering, Princeton University, Princeton, NJ 08544, USA}

\begin{document}

\ninept
\maketitle

\begin{abstract}
The paper deals with minimax optimal statistical tests for two composite hypotheses, where each hypothesis is defined by a non-parametric uncertainty set of feasible distributions. It is shown that for every pair of uncertainty sets of the $f$-divergence ball type, a pair of uncertainty sets of the density band type can be constructed, which is equivalent in the sense that it admits the same pair of least favorable distributions. This result implies that robust tests under $f$-divergence ball uncertainty, which are typically only minimax optimal for the single sample case, are also fixed sample size minimax optimal with respect to the equivalent density band uncertainty sets.
\end{abstract}

%%%%%%%%%%%%%%%%%%%%%%%%%%%%%%%%%%%%%%%%%%%%%%%%%%%%%%%%%%%%%%%%%%%%%%%
\section{Introduction}
\label{sec:intro}
%%%%%%%%%%%%%%%%%%%%%%%%%%%%%%%%%%%%%%%%%%%%%%%%%%%%%%%%%%%%%%%%%%%%%%%

A statistical test for two composite hypotheses is called minimax optimal if it minimizes the maximum risk over the two corresponding sets of feasible distributions. In the context of robust statistics, these sets are referred to as \emph{uncertainty sets}. In contrast to adaptive procedures \cite{Zeitouni1992_glrt}, the minimax approach provides strict guarantees on the error probabilities for all feasible distributions. Moreover, minimax tests are often easy to implement since they typically reduce to an optimal test for two simple hypotheses, where each hypothesis is represented by a \emph{least favorable distribution}.

A common way of specifying uncertainty sets is via a neighborhood around a nominal distribution, which represents an ideal system state or model \cite{Kassam1981_robustness_survey}. In many works on robust detection, the use of $f$-divergence balls has been proposed as a useful and versatile model to construct such neighborhoods \cite{McKellips_binary_input, McKellips1998_maximin, Levy2009_entropy_tolerance, Gul2013_modelling_errors, Gul2014_Hellinger_distance, Gul2015_composite_distances, Gul2016_alpha_divergence, Gul2017_minimax_robust}. In contrast to outlier models, such as $\varepsilon$-contamination \cite{Huber1965_robust_PRT}, $f$-divergence balls do not allow for arbitrarily large deviations from the nominals and, therefore, have been argued to better represent scenarios where the shape of a distribution is subject to uncertainty, but there are no gross outliers in the data \cite{Levy2009_entropy_tolerance}.

In order to present the result in this paper, the concept of single sample and fixed sample size tests needs to be introduced. A single sample test is based on the observation of a single, possibly vector-valued, random variable $X_1$. Consequently, the uncertainty sets are defined in terms of all possible \emph{joint} distributions of the elements of $X$. Such an uncertainty model is suitable in some cases, but more often the observations are obtained by repeatedly performing independent random experiments so that the test is based on a sequence of independent random variables $X_1, \ldots, X_N$, $N > 1$. By definition, this independence constraint cannot be incorporated into a single sample test. Hence, tests whose observations are realizations of multiple independent random variables need to be considered separately. In order to highlight the difference to tests whose sample size is random \cite{Wald1947_sequential_analysis}, they are referred to as fixed sample size tests in what follows.

For most commonly used uncertainty models---including the density band model, which will be discussed in detail later on---it can be shown that a single sample minimax optimal test is also fixed sample size minimax optimal. More precisely, the least favorable distributions for $X_1$ in the single sample case are also least favorable for all $X_1, \ldots, X_N$ in the fixed sample size case. However, in general, this does not hold true for uncertainty sets of the $f$-divergence ball type, where the fixed sample size minimax optimal solution is typically intractable. However, a commonly applied heuristic is to use the single sample least favorable distributions for the fixed sample size case anyhow, regardless of the fact that this extension does not hold in theory; compare \cite[Sec.~V.A]{Gul2017_minimax_robust}. Tests constructed this way are referred to as single sample minimax optimal tests with repeated observations. Evidently, such tests are no longer minimax optimal. However, in practice, it can be observed that they are still robust, meaning that they meet the specified error probabilities for most if not all feasible distributions.

In this paper, the favorable robustness properties of single sample minimax optimal tests with repeated observations are explained in a rigorous manner. It is shown that they are indeed fixed sample size minimax optimal, but for a density band uncertainty model instead of the original $f$-divergence ball model. That is, single-sample minimax optimal tests can be applied to repeated observations without sacrificing minimax optimality, if one is willing to accept a change in the uncertainty model. This result is proved by showing that for every $f$-divergence ball model, there exists and equivalent density band model that admits the same single sample minimax optimal solution. However, for the density band model, this solution is known to be fixed sample size optimal as well so that it automatically extends to the case of repeated observations.

The paper is organized as follows: a brief review of minimax optimal detection is given in Section~\ref{sec:minimax_optimal_detection}. The two uncertainty models of interest, i.e., $f$-divergence balls and density bands, are introduced in Section~\ref{sec:uncertainty_sets}. The main result is stated and proved in Section~\ref{sec:main_result}, followed by a brief discussion in Section~\ref{sec:discussion}. An illustrative example is shown in Section~\ref{sec:example}, which also concludes the paper.

%%%%%%%%%%%%%%%%%%%%%%%%%%%%%%%%%%%%%%%%%%%%%%%%%%%%%%%%%%%%%%%%%%%%%%%
\section{Minimax Optimal Detection}
\label{sec:minimax_optimal_detection}
%%%%%%%%%%%%%%%%%%%%%%%%%%%%%%%%%%%%%%%%%%%%%%%%%%%%%%%%%%%%%%%%%%%%%%%

The single sample case is considered first. Let $(\mathcal{X},\mathcal{F})$ be a measurable space and let $X_1$ be a $(\mathcal{X},\mathcal{F})$-valued random variable that is distributed according to a probability measure (distribution) $P$. Throughout the paper it is assumed that all distributions on $(\mathcal{X},\mathcal{F})$ have a continuous density function with respect to some $\sigma$-finite reference measure $\mu$. The set of all distributions on $(\mathcal{X},\mathcal{F})$ that admit this property is denoted by $\mathcal{M}_\mu$. %\pagebreak

The goal of a simple, non-robust binary hypothesis test is to decide between the two hypotheses
\begin{align*}
  \mathcal{H}_0\colon \; P = P_0, \qquad
  \mathcal{H}_1\colon \; P = P_1,
\end{align*}
where $P_0, P_1 \in \mathcal{M}_\mu$ are two given distributions. The test is defined by a decision $d \in \{0,1\}$ and a, possibly randomized, decision rule $\delta\colon \mathcal{X} \to [0,1]$, where $\delta(x)$ denotes the conditional probability to decide for $\mathcal{H}_1$, given the observation $X_1 = x$. The set of all decision rules is denoted by $\Delta$. The type I and type II error probabilities are given by
\begin{align*}
  P_0[d=1] &= E_{P_0}[\,\delta(X) \,], \\
  P_1[d=0] &= E_{P_1}[1-\delta(X)].
\end{align*}
The optimal decision rule $\delta^*$ for the simple binary hypothesis test is a threshold comparison of the likelihood ratio, i.e.,
\begin{equation*}
  \delta^*(x) = \begin{cases}
                  1,      & l(x) > \lambda \\
                  \kappa, & l(x) = \lambda \\
                  0,      & l(x) < \lambda
                \end{cases},
\end{equation*}
where $\lambda > 0$ is the threshold value, $\kappa \in [0,1]$ can be chosen arbitrarily, and $l(x)$ denotes the likelihood ratio
\begin{equation*}
  l(x) = \frac{p_1(x)}{p_0(x)}.
\end{equation*}
The likelihood ratio test is optimal in a very general sense \cite{Christensen2005_Fisher-Neyman-Pearson-Bayes}. In particular, it minimizes the weighted sum error probability, i.e., it solves
\begin{equation}
  \label{eq:decision_rule_simple}
  \min_{\delta \in \Delta} \; E_{P_1}[\,\delta(X) \,] + \lambda \, E_{P_0}[1-\delta(X)].
\end{equation}

In robust detection, the distribution under each hypothesis is assumed not to be known exactly. The distributional uncertainty is modeled by two disjoint sets $\mathcal{P}_0,\mathcal{P}_1 \subset \mathcal{M}_\mu$ so that the hypotheses become
\begin{align*}
  \mathcal{H}_0\colon \; P \in \mathcal{P}_0, \qquad
  \mathcal{H}_1\colon \; P \in \mathcal{P}_1.
\end{align*}
The minimax problem corresponding to \eqref{eq:decision_rule_simple} is thus given by
\begin{equation}
  \label{eq:minimax_problem}
  \min_{\delta \in \Delta} \; \max_{\substack{H_0 \in \mathcal{P}_0 \\ H_1 \in \mathcal{P}_1}} \; E_{H_1}[\, \delta(X) \,] + \lambda \, E_{H_0}[1-\delta(X)].
\end{equation}
Problem \eqref{eq:minimax_problem} is central to robust detection. By definition, its solution is minimax optimal with respect to the weighted sum of error probabilities, but it can also be shown to be minimax optimal in the sense of Neyman--Pearson and the Baysian sense. This property is fixed in the following definition.
\begin{definition}
  A triplet $(\delta^*,Q_0,Q_1)$ that solves \eqref{eq:minimax_problem} for a given $\lambda > 0$ is called \emph{single sample minimax optimal} with respect to the threshold $\lambda$ and the uncertainty sets $\mathcal{P}_0$, $\mathcal{P}_1$. This is written as
  \begin{equation*}
    (\delta^*,Q_0,Q_1) \in \{\mathcal{P}_0,\mathcal{P}_1\}_{\lambda}^*.
  \end{equation*}
\end{definition}
In \cite{Huber1965_robust_PRT} and \cite{Fauss2016_old_bands} it is shown that if the solution of \eqref{eq:minimax_problem} is independent of the threshold $\lambda$, it is also minimax optimal for fixed sample size tests with arbitrary thresholds and arbitrary sample sizes. This property is fixed in the next definition.
\begin{definition}
  A triplet $(\delta^*,Q_0,Q_1)$ that jointly solves \eqref{eq:minimax_problem} for all $\lambda > 0$ is called \emph{fixed sample size minimax optimal} with respect to $\mathcal{P}_0$, $\mathcal{P}_1$. This is written as
  \begin{equation*}
    (\delta^*,Q_0,Q_1) \in \{\mathcal{P}_0,\mathcal{P}_1\}^*.
  \end{equation*}
\end{definition}

%%%%%%%%%%%%%%%%%%%%%%%%%%%%%%%%%%%%%%%%%%%%%%%%%%%%%%%%%%%%%%%%%%%%%%%
\section{Uncertainty Sets}
\label{sec:uncertainty_sets}
%%%%%%%%%%%%%%%%%%%%%%%%%%%%%%%%%%%%%%%%%%%%%%%%%%%%%%%%%%%%%%%%%%%%%%%

Two types of uncertainty sets are introduced in this section. The first one is the $f$-divergence ball model, which specifies uncertainty sets in terms of a maximum feasible distance from a nominal distribution and allows the use of arbitrary $f$-divergences to define this distance. Formally, $f$-divergence ball uncertainty sets are of the form
\begin{equation}
  \label{eq:f-divergence_uncertainty}
  \mathcal{P}_f(P,\varepsilon) = \{ H \in \mathcal{M}_{\mu} : D_f(H \Vert P) \leq \varepsilon \},
\end{equation}
where $P$ denotes the nominal distribution and $D_f$ denotes the $f$-divergence induced by the function $f$, i.e.,
\begin{align*}
  D_f(H \Vert P) &= \int_{\mathcal{X}} f\biggl( \frac{\mathrm{d} H}{\mathrm{d} P}(x) \biggr) \, \mathrm{d} P(x) \\
  &= \int_{\mathcal{X}} f\biggl( \frac{h(x)}{p(x)} \biggr) p(x) \, \mathrm{d} \mu(x),
\end{align*}
where $f \colon \mathbb{R}_{\geq 0} \to \mathbb{R}$ is convex and satisfies $f(1) = 0$. The definition of the $f$-divergence ball in terms of $D_f(H \Vert P)$ instead of $D_f(P \Vert H)$ is arbitrary since for every feasible function $f$ it holds that $D_f(P \Vert H) = D_{\tilde{f}}(H \Vert P)$, with $\tilde{f}(x) = f(\tfrac{1}{x})x$.

Owing to the mild constraints on $f$, uncertainty sets of the form \eqref{eq:f-divergence_uncertainty} offer a great amount of flexibility and have attracted increased attention in recent years. In \cite{Levy2009_entropy_tolerance} and \cite{Gul2017_minimax_robust}, minimax optimal tests based on the Kullback--Leibler divergence were derived under varying assumptions. Minimax optimal tests have also been derived for the squared Hellinger distance \cite{Gul2013_modelling_errors, Gul2014_Hellinger_distance}, the total variation distance \cite{Gul2015_composite_distances}, and $\alpha$-divergences \cite{Gul2016_alpha_divergence, Gul2017_minimax_robust}. However, a disadvantage of robust tests with $f$-divergence ball uncertainty is that no fixed sample size minimax optimal solution is guaranteed to exist. In fact, most of the works cited above only consider the single sample case.

The second type of uncertainty model is the density band model. In a robust detection context, it was first proposed by Kassam \cite{Kassam1981} and covers sets of the form
\begin{equation}
  \label{eq:density_band_uncertainty}
  \mathcal{P}_\text{b}(P',P'') = \{ H \in \mathcal{M}_{\mu} : p' \leq h \leq p'' \},
\end{equation}
where $P',P''$ are nonnegative measures on $(\mathcal{X},\mathcal{F})$ that admit densities $p',p''$ with respect to $\mu$ and satisfy
\begin{equation*}
  P'(\mathcal{X}) \leq 1, \quad  P''(\mathcal{X}) \geq 1, \quad \text{and} \quad 0 \leq p' \leq p''.
\end{equation*}
In words, the density band model restricts the true density to lie within a band specified by $p'$ and $p''$. Similar to the choice of $f$ in the $f$-divergence ball model, the choice of $P'$ and $P''$ allows for \emph{a priori} knowledge about the type of contamination to be incorporated into the model. Therefore, although it is still based on the concept of outliers, the band model can capture a much larger variety of contamination types and mismatches than the standard $\varepsilon$-contamination model. Another useful property of the density band model is that for every pair of uncertainty sets of the form \eqref{eq:density_band_uncertainty}, a fixed sample size minimax optimal test is guaranteed to exist. Moreover, the corresponding least favorable densities can be calculated in a generic manner using a simple and efficient algorithm. See \cite{Fauss2016_old_bands} for a more detailed discussion of the band model and the calculation of its least favorable densities.

In the next section, it is shown in how far the two uncertainty models can be considered equivalent in the single sample case and how the density band model can be used to construct fixed sample size minimax optimal tests from tests that are only single sample size minimax optimal under $f$-divergence ball uncertainty.

%%%%%%%%%%%%%%%%%%%%%%%%%%%%%%%%%%%%%%%%%%%%%%%%%%%%%%%%%%%%%%%%%%%%%%%
\section{Main Result}
\label{sec:main_result}
%%%%%%%%%%%%%%%%%%%%%%%%%%%%%%%%%%%%%%%%%%%%%%%%%%%%%%%%%%%%%%%%%%%%%%%

In this section, the main result of the paper is stated and proved. A more detailed discussion is deferred to Section~\ref{sec:discussion}.
\begin{theorem*}
  Let $\mathcal{P}_{f_0}(P_0,\varepsilon_0)$ and $\mathcal{P}_{f_1}(P_1,\varepsilon_1)$ be two uncertainty sets of the form \eqref{eq:f-divergence_uncertainty}. If it holds that
  \begin{equation*}
    (\delta^*, Q_0, Q_1) \in \{\mathcal{P}_{f_0}(P_0,\varepsilon_0) \,,\, \mathcal{P}_{f_1}(P_1,\varepsilon_1)\}_{\lambda}^*,
  \end{equation*}
  then there exist nonnegative scalars $a_0 \leq b_0$ and $a_1 \leq b_1$ such that
  \begin{equation*}
    (\delta^*, Q_0, Q_1) \in \{\mathcal{P}_\text{b}(a_0 P_0,b_0 P_0) \,,\, \mathcal{P}_\text{b}(a_1 P_1,b_1 P_1)\}^*.
  \end{equation*}
\end{theorem*}

In words, the theorem states that if $(\delta^*,Q_0,Q_1)$ is single sample minimax optimal for an $f$-divergence ball uncertainty model, a density band model can be constructed from scaled versions of the nominal densities such that $(\delta^*,Q_0,Q_1)$ is fixed sample size minimax optimal with respect to this band model. A proof is detailed below.

\begin{proof*}
Let $P_0$, $P_1$, $f_0$, $f_1$, $\varepsilon_0$, $\varepsilon_1$, and $\lambda$ be given. Rewriting the optimization problem \eqref{eq:minimax_problem} in terms of the densities and with explicit constraints yields
\begin{gather}
  \max_{\substack{h_0 > 0 \\ h_1 > 0}} \; \min_{\delta \in \Delta} \; \int_{\mathcal{X}} h_1 \, \delta + \lambda h_0 (1-\delta) \, \mathrm{d}\mu
  \quad \text{s.t.} \label{eq:minimax_objective} \\
  \int_{\mathcal{X}} f_0\biggl( \frac{h_0}{p_0} \biggr) p_0 \, \mathrm{d}\mu \leq \varepsilon_0, \quad
  \int_{\mathcal{X}} f_1\biggl( \frac{h_1}{p_1} \biggr) p_1 \, \mathrm{d}\mu \leq \varepsilon_1 \label{eq:f-divergence_constraints} \\
  \int_{\mathcal{X}} h_0 \, \mathrm{d}\mu = 1, \quad
  \int_{\mathcal{X}} h_1 \, \mathrm{d}\mu = 1. \label{eq:density_constraints}
\end{gather}
By assumption, $(\delta^*,Q_0,Q_1)$ solves this minimax problem, which implies that $(\delta^*,Q_0,Q_1)$ is a saddle point of \eqref{eq:minimax_objective}. This, in turn, implies that $(\delta^*,Q_0,Q_1)$ satisfies the corresponding Karush--Kuhn--Tucker conditions, which are first order necessary conditions for optimality \cite{Guignard1969}. In particular, stationarity of the saddle point solution implies that scalars $\eta_0,\eta_1$ and nonnegative scalars $\nu_0,\nu_1$ exists such that
\begin{align}
  \lambda(1-\delta^*) &= \nu_0 f_0'\biggl( \frac{q_0}{p_0} \biggr) - \eta_0, \label{eq:stationarity_0}\\
  \delta^* &= \nu_1 f_1'\biggl( \frac{q_1}{p_1} \biggr) - \eta_1, \label{eq:stationarity_1} \\
  \delta^* &= \begin{cases}
              1,                 & q_1 > \lambda q_0 \\
              \kappa \in (0,1),  & q_1 = \lambda q_0 \\
              0,                 & q_1 < \lambda q_0
            \end{cases}, \label{eq:stationarity_delta}
\end{align}
where $f_0'$, $f_1'$ denote the (sub)derivatives of $f_0$ and $f_1$, $\eta_0$, $\eta_1$ denote the Lagrange multipliers corresponding to the constraints \eqref{eq:density_constraints}, and $\nu_0$, $\nu_1$ denote the Lagrange multipliers corresponding to the constraints \eqref{eq:f-divergence_constraints}. Since $f_0$ and $f_1$ are convex, their (sub)derivatives are nondecreasing. Moreover, the inverse functions $g_0$ and $g_1$, which are implicitly defined by
\begin{equation*}
  g_0(f_0'(x)) = x \quad \text{and} \quad g_1(f_1'(x)) = x \quad \forall x \in \mathbb{R}_{\geq 0},
\end{equation*}
exist and are nondecreasing as well. Solving \eqref{eq:stationarity_0} and \eqref{eq:stationarity_1} for $q_0$ and $q_1$ yields
\begin{align}
  q_0 &= g_0\biggl(\frac{\lambda(1-\delta^*) + \eta_0}{\nu_0}\biggr) p_0, \label{eq:q0_g0} \\
  q_1 &= g_1\biggl(\frac{\delta^* + \eta_1}{\nu_1}\biggr) p_1. \label{eq:q1_g1}
\end{align}
Combining \eqref{eq:stationarity_delta}, \eqref{eq:q0_g0} and \eqref{eq:q1_g1}, it follows that the least favorable densities are of the form
\begin{equation}
  q_0 = \begin{cases}
          b_0 \, p_0,               & \delta^* = 0 \\
          \frac{1}{\lambda} \, q_1, & \delta^* \in (0,1) \\
          a_0 \, p_0,               & \delta^* = 1
        \end{cases}
  \label{eq:lfds_piecewise_0}
\end{equation}
and
\begin{equation}
  \label{eq:lfds_piecewise_1}
  q_1 = \begin{cases}
          a_1 \, p_1,      & \delta^* = 0 \\
          \lambda \, q_0,  & \delta^* \in (0,1) \\
          b_1 \, p_1,      & \delta^* = 1
        \end{cases},
\end{equation}
where
\begin{equation}
  \label{eq:coefficients_a}
  a_0 = g_0\left(\frac{\eta_0}{\nu_0}\right) \leq b_0 = g_0\left(\frac{\lambda + \eta_0}{\nu_0} \right)
\end{equation}
and
\begin{equation}
  \label{eq:coefficients_b}
  a_1 = g_1\left(\frac{\eta_1}{\nu_1}\right) \leq b_1 = g_1\left(\frac{1 + \eta_1}{\nu_1}\right).
\end{equation}
Note that since $q_0$ and $q_1$ are valid densities, $a_0, b_0$ and $a_1, b_1$ are nonnegative. The next step of the proof is to show that the least favorable densities in \eqref{eq:lfds_piecewise_0} and \eqref{eq:lfds_piecewise_1} can be written as
\begin{align}
  q_0 &= \min \bigl\{ b_0 p_0 \,,\, \max\bigl\{ \tfrac{1}{\lambda} q_1 \,,\, a_0 p_0 \bigr\} \bigr\}, \label{eq:lfd0} \\
  q_1 &= \min \bigl\{ b_1 p_1 \,,\, \max\bigl\{ \lambda q_0 \,,\, a_1 p_1 \bigr\} \bigr\}. \label{eq:lfd1}
\end{align}
Only \eqref{eq:lfd0} is shown here since the proof for \eqref{eq:lfd1} can be given analogously. From \eqref{eq:stationarity_delta} and \eqref{eq:lfds_piecewise_0} it follows that on $\{ x \in \mathcal{X} : \delta(x) = 0\}$
\begin{align}
  q_1 < \lambda q_0 \quad &\Rightarrow \quad \frac{1}{\lambda} q_1 < q_0 = b_0 p_0, \label{eq:lfd0_d0} \\
  \intertext{on $\{ x \in \mathcal{X} : \delta(x) \in (0,1)\}$}
  q_1 = \lambda q_0 \quad &\Rightarrow \quad \frac{1}{\lambda} q_1 = q_0, \label{eq:lfd0_d01} \\
  \intertext{and on $\{ x \in \mathcal{X} : \delta(x) = 1\}$}
  q_1 > \lambda q_0 \quad &\Rightarrow \quad \frac{1}{\lambda} q_1 > q_0 = a_0 p_0. \label{eq:lfd0_d1}
\end{align}
Combining \eqref{eq:lfd0_d0}, \eqref{eq:lfd0_d01}, and \eqref{eq:lfd0_d1} yields \eqref{eq:lfd0}. From \cite[Theorem 4]{Fauss2016_old_bands}, it follows immediately that \eqref{eq:lfd0} and \eqref{eq:lfd1} are fixed sample size least favorable for a density band model with bounds
\begin{align*}
  p_0' &= a_0 p_0 \leq b_0 p_0 = p_0'', \\
  p_1' &= a_1 p_1 \leq b_1 p_1 = p_1''.
\end{align*}
\end{proof*}

%%%%%%%%%%%%%%%%%%%%%%%%%%%%%%%%%%%%%%%%%%%%%%%%%%%%%%%%%%%%%%%%%%%%%%%
\section{Discussion}
\label{sec:discussion}
%%%%%%%%%%%%%%%%%%%%%%%%%%%%%%%%%%%%%%%%%%%%%%%%%%%%%%%%%%%%%%%%%%%%%%%

The result presented in the previous section states that every single sample minimax optimal test under $f$-divergence ball uncertainty is fixed sample size minimax optimal under the equivalent density band uncertainty. This not only makes it possible to use single sample results for fixed sample size tests without sacrificing minimax optimality, but also to specify the exact sets of distributions for which the minimax property holds. In this sense, the theorem lifts the $f$-divergence ball model to the same level of usefulness as the classic outlier models, whose single sample results automatically carry over to the fixed sample case. In addition to this generalization, the fact that for every $f$-divergence ball model an equivalent density band model can be constructed offers some deeper insights and also suggests an alternative approach to the design of robust tests under $f$-divergence ball uncertainty.

One useful aspect of the equivalent band model is that it simplifies comparing the amount and type of uncertainty that is allowed for by different $f$-divergence ball models. Such comparisons are non-trivial since the $\varepsilon$-tolerances in \eqref{eq:f-divergence_uncertainty} do not directly translate to contamination ratios and might be of different scales altogether. While, for example, the Kullback--Leibler divergence can take on any nonnegative value, the Hellinger distance is bounded between zero and one. In such cases, one cannot simply compare the $\varepsilon$-tolerances in order to compare the maximum amount of uncertainty in the distributions. The corresponding band model, however, offers a way to make such comparisons possible. The lower bounds on the density functions, $a_0 p_0$ and $a_1 p_1$, determine how much probability mass the nominal distributions contribute at least, namely $a_0$ and $a_1$. Consequently, the outliers can at most contribute the remaining probability masses $1-a_0$ and $1-a_1$, which can hence be interpreted as contamination ratios. The larger they are, the more uncertainty a model allows. On the other hand, the upper bounds, $b_0 p_0$ and $b_1 p_1$, offer an insight into what type of uncertainty is allowed. For $b_0, b_1 \gg 1$ the contamination is almost unconstrained, which corresponds to gross outliers. For $b_0,b_1 \approx 1$, the outlier distributions are close to the nominals, which corresponds to more subtle model mismatches. In general, the outlier distributions under $\mathcal{H}_i$, $i \in \{0,1\}$, are constrained to lie within the set
\begin{equation*}
  \left\{ H \in \mathcal{M}_\mu : h \leq \frac{b_i-a_i}{1-a_i} p_i \right\}.
\end{equation*}
This interpretation of a density band model as a constrained $\varepsilon$-contamination model often offers a useful intuition for the amount and type of contamination that cannot be obtained by inspection of the $f$-divergence ball model.

There are several ways to determine the coefficients $a_0$, $b_0$ and $a_1$, $b_1$ in practice. If expressions for the least favorable distributions can already be found in the literature, the coefficients can be determined by a simple comparison with the expressions in \eqref{eq:lfds_piecewise_0} and \eqref{eq:lfds_piecewise_1}. If the least favorable densities are unknown, the KKT conditions \eqref{eq:stationarity_0}--\eqref{eq:stationarity_delta} can be solved for $\nu_0$ ,$\nu_1$ and $\eta_0$, $\eta_1$. The scaling coefficients can then be calculated according to \eqref{eq:coefficients_a} and \eqref{eq:coefficients_b}. In practice, however, this approach might be prohibitively complex.

An alternative to solving the KKT conditions for $\nu_0,\nu_1$ and $\eta_0,\eta_1$ is to solve them directly for $a_0,b_0$ and $a_1,b_1$. From the result in the previous section, it follows that the least favorable densities are of the form \eqref{eq:lfd0} and \eqref{eq:lfd1}. For given scaling coefficients $a_0,b_0$ and $a_1,b_1$, these equations can be solved for $q_0$ and $q_1$ by finding a threshold $\lambda$ so that the right hand sides of \eqref{eq:lfd0} and \eqref{eq:lfd1} are valid densities, i.e., they integrate to one. Finally, an outer search over $a_0,b_0$ and $a_1,b_1$ can be performed such that the primal constraints are fulfilled, i.e.,
\begin{equation*}
  \int_{\mathcal{X}} f_0 \biggl( \frac{q_0}{p_0} \biggr) p_0 \, \mathrm{d} \mu = \varepsilon_0
  \quad \text{and} \quad
  \int_{\mathcal{X}} f_1 \biggl( \frac{q_1}{p_1} \biggr) p_1 \, \mathrm{d} \mu = \varepsilon_1.
\end{equation*}
This approach can be expected to be less efficient than a solution that exploits properties of a given function $f$, but is applicable in general and does not require prior analysis of the problem.

Yet another option to determine the scaling coefficients is to directly solve the primal problem \eqref{eq:minimax_objective} using a suitable convex optimization algorithm. Even if the latter does not calculate the dual variables explicitly, $a_0,b_0$ and $a_1,b_1$ can be obtained from the ratio of the least favorable and the nominal densities
\begin{equation}
  \label{eq:lfd_ratio}
  \begin{aligned}
    \frac{q_0}{p_0} &= \min \left\{ b_0 \,,\, \max \left\{ \frac{1}{\lambda} \frac{q_1}{p_0} \,,\, a_0 \right\} \right\}, \\
    \frac{q_1}{p_1} &= \min \left\{ b_1 \,,\, \max \left\{ \lambda \frac{q_0}{p_1} \,,\, a_1 \right\} \right\}.
  \end{aligned}
\end{equation}
By inspection of \eqref{eq:lfd_ratio}, the scaling coefficients can be identified form the regions where the likelihood ratio is constant.

The result in Section~\ref{sec:main_result} also motivates further research into how the two uncertainty models are related. Given equivalent $f$-divergence balls and density bands, does one uncertainty set contain the other, i.e., is one model a relaxation of the other one? Does every band model whose bounds are constructed by scaling a nominal density admit an equivalent $f$-divergence ball model? If so, how can the corresponding function $f$ be determined? Another question that might be asked is, whether similar equivalences exist for other types of uncertainty models as well and in how far there is a hierarchy between these models, i.e., whether the set of all possible least favorable distributions under one model is a sub- or super-set of all possible least favorable distributions under another model.

%%%%%%%%%%%%%%%%%%%%%%%%%%%%%%%%%%%%%%%%%%%%%%%%%%%%%%%%%%%%%%%%%%%%%%%
\section{Example}
\label{sec:example}
%%%%%%%%%%%%%%%%%%%%%%%%%%%%%%%%%%%%%%%%%%%%%%%%%%%%%%%%%%%%%%%%%%%%%%%

\begin{figure}[t]
  \centering
  \includegraphics{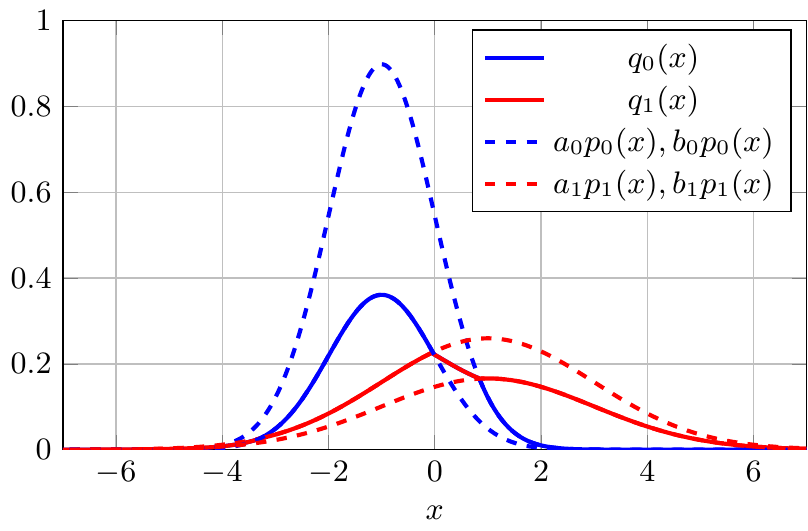}
  \caption{Least favorable densities and equivalent density bands for uncertainty sets  $\mathcal{P}_{x\log}\bigl(\mathcal{N}(-1,1),0.03\bigr)$ and $\mathcal{P}_{x\log}\bigl(\mathcal{N}(1,2),0.02\bigr)$.}
  \label{fig:bands}
\end{figure}

In order to highlight the connection to existing results, we consider the example from \cite[Sec.~VI.A]{Gul2017_minimax_robust}, where the nominal distributions under $\mathcal{H}_0$ and $\mathcal{H}_1$ are chosen as $P_0 = \mathcal{N}(-1,1)$ and $P_1 = \mathcal{N}(1,2)$, respectively, and $\mathcal{N}(m,\sigma^2)$ denotes a Gaussian distribution with mean $m$ and variance $\sigma^2$. The Kullback--Leibler divergence is used as a distance measure, i.e., $f_0(x) = f_1(x) = x\log(x)$, the tolerances are chosen as $\varepsilon_0 = 0.03$, $\varepsilon_1 = 0.02$, and the likelihood ratio threshold is set to $\lambda = 1$. Using Theorem~2 in \cite{Gul2017_minimax_robust}, the least favorable densities for this model can be calculated efficiently by solving two integral equations. The coefficients for the corresponding band model can be identified by comparing (6) in \cite{Gul2017_minimax_robust} to \eqref{eq:lfds_piecewise_0} and \eqref{eq:lfds_piecewise_1} in this paper. For the numerical values given above, they calculate to $a_0 \approx 0.9047$, $b_0 \approx 2.2519$, $a_1 \approx 0.8319$, and $b_1 \approx 1.3009$. The least favorable densities and the equivalent density bands are depicted in Fig.~\ref{fig:bands}. Interestingly, the difference in the radii of the $f$-divergence balls is not reflected in the contamination ratio, which is smaller under $\mathcal{H}_0$ ($\approx 10\%$) than under $\mathcal{H}_1$ ($\approx 17\%$). However, the band under $\mathcal{H}_0$ is wider, meaning that it allows for larger deviations from the nominal distribution. Under $\mathcal{H}_1$, the contamination ratio is higher, but the outlier distribution is much more restricted.

This example illustrates how existing results on robust tests under $f$-divergence ball uncertainty can be used to construct fixed sample size minimax optimal tests for the equivalent density band uncertainty sets and how the latter provide additional insight into the amount and type of contamination induced by the uncertainty model.

% References should be produced using the bibtex program from suitable
% BiBTeX files (here: strings, refs, manuals). The IEEEbib.bst bibliography
% style file from IEEE produces unsorted bibliography list.
% -------------------------------------------------------------------------
\bibliographystyle{IEEEbib}
\bibliography{references}
\end{document}